\documentclass[11pt]{article}
\usepackage{amssymb, amsmath,amsthm,amsfonts,a4,bezier}
\usepackage{hyperref}
\usepackage{amscd}
\usepackage{graphpap}
\usepackage[pdftex]{color,graphicx,xcolor}
\usepackage{color,graphicx,enumerate}
\usepackage{mathpazo}
\usepackage{setspace}
\usepackage{comment}
\usepackage{mathtools}
\usepackage{relsize}

\DeclareGraphicsExtensions{.jpg,.pdf,.eps,.png}
\newtheorem{theo}{Theorem}

\newtheorem{Rem}{Remark}

\newtheorem{Conjecture}{Conjecture}
\renewcommand{\o}{{\omega}}

\renewcommand{\u}{{\mathbf {u}}}
\newcommand{\Z}{{\mathbb {Z}}}

\newcommand{\E}{{\mathbb {E}}}
\newcommand{\F}{{\mathcal {F}}}

\newcommand{\PP}{{\mathbb{P}}}
\newcommand{\N}{{\mathbb {N}}}
\begin{document}
\thispagestyle{empty}
\baselineskip=28pt
\vskip 5mm
\begin{center} {\LARGE{\bf A Note on the Dimensional Crossover Critical Exponent}}
\end{center}

\baselineskip=12pt
\vskip 10mm

\begin{center}\large
Pablo A. Gomes\footnote{Departamento de Matem\'{a}tica, Universidade Federal de Minas Gerais, Brazil, \href{mailto:rsanchis@mat.ufmg.br}{pabloag@ufmg.br}; https://orcid.org/0000-0002-4472-3628} , R\'{e}my Sanchis\footnote{Departamento de Matem\'{a}tica, Universidade Federal de Minas Gerais, Brazil, \href{mailto:rsanchis@mat.ufmg.br}{rsanchis@mat.ufmg.br;} https://orcid.org/0000-0002-3476-7459
} 
and Roger W. C. Silva 
 \footnote{Corresponding Author - Departamento de Estat\'{\i}stica,
 Universidade Federal de Minas Gerais, Brazil, \href{mailto:rogerwcs@ufmg.br} {rogerwcs@ufmg.br;} https://orcid.org/0000-0002-6365-9211
}
\\ 


\end{center}

\begin{abstract} We consider independent anisotropic bond percolation on $\Z^d\times \Z^s$ where edges parallel to $\Z^d$ are open with probability $p<p_c(\Z^d)$ and edges parallel to $\Z^s$ are open with probability $q$, independently of all others. We prove that percolation occurs for $q\geq 8d^2(p_c(\Z^d)-p)$. This fact implies that the so-called {\it Dimensional Crossover}  critical exponent, if it exists, is greater or equal than 1. In particular, using known results, we conclude the proof that, for $d\geq 11$, the crossover  critical exponent exists and equals 1.\\


\end{abstract}

\onehalfspacing
\section{Introduction and Results}
\subsection{Background}\label{background}

\hspace{0.5cm}In this note we consider anisotropic bond percolation on the graph $(\Z^{d+s},E(\Z^{d+s}))$, where $E(\Z^{d+s})$ is the set of edges between nearest neighbors of $\Z^{d+s}$. We simplify notation and denote this graph by $\Z^{d+s}=\Z^d\times \Z^s$. An edge of $\Z^{d+s}$ is called a $\Z^d$-edge (respectively a $\Z^s$-edge) if it joins two vertices which differ only in their $\Z^d$ (respectively $\Z^s$) component. Formally a $\Z^d$-edge is of the form $\langle(u_d,u_s),(u'_d,u_s)\rangle$, where $u_d$ and $u_s$ should be understood as the components of $u\in\Z^{d+s}$. Probability is introduced as follows:  given two parameters $p,q\in[0,1]$, we declare each $\Z^d$-edge open with probability $p$ and each $\Z^s$-edge open with probability $q$, independently of all others. This model is described by the probability space $(\Omega, \F, \mathbb{P}_{p,q})$ where
$\Omega = \{0,1\}^{E(\Z^{d+s})}$, $\F$ is the $\sigma$-algebra generated by the cylinder sets in $\Omega$ and $\mathbb{P}_{p,q} =\prod_{e\in E}\mu(e)$, where $\mu(e)$ is Bernoulli measure with parameter $p$ or $q$ according to $e$ been a $\Z^d$- edge or a $\Z^s$- edge respectively.

Given two vertices $u,v\in\Z^{d+s}$, we say that $u$ and $v$ are connected in the configuration $\o$  if there exists an open path in $\Z^{d+s}$ starting in $u$ and ending in $v$. The event where $v$ and $u$ are connected is denoted by $\{\o\in\Omega:v\leftrightarrow u\mbox{ in $\o$}\}$ and we write $C(\o)=\{u\in \Z^{d+s}:u\leftrightarrow 0\mbox{ in $\o$}\}$ for the open cluster containing the origin. We denote by $\theta(p,q)=\mathbb{P}_{p,q}(\o\in\Omega:|C(\o)|=\infty)$ the main macroscopic function in percolation theory and denote the mean size of the open cluster by $\chi(p,q)=\E_{p,q}(|C(\o))|)$. Whenever necessary we shall write  $\chi_p(d)$  and $p_c(d)$  for the expected cluster size and critical threshold on $\Z^d$  with a single parameter $p\in(0,1)$. For a comprehensive background in percolation theory we refer the reader to \cite{GRI}.

It is easy to see, by a standard coupling argument,  that $\theta(p,q)$ is a monotone non-decreasing function of the parameters $p$ and $q$. This enables us to define the function $q_c:[0,1]\rightarrow [0,1]$, where
\begin{equation}\label{qcrit}q_c(p)=\sup\{q:\theta(p,q)=0\}.
\end{equation}
The function $q_c(p)$ is continuous and strictly decreasing (see \cite{CLS}) and we are interested in understanding its  behavior  as $p\uparrow p_c(d)$? 

This problem arises in the physics literature as the dimensional crossover problem. The term crossover relates to the study of percolative systems on  $(d+s)$-dimensional lattices, where the $d$-dimensional parameter $p$ is close to $p_c(d)$ from below and the $s$-dimensional parameter $q$ is small. Similar anisotropic ferromagnectic models have also been considered in the mathematical physics literature, see the works \cite{FMMPV}, \cite{FMMPV2} and  \cite{MPS}.

\subsection{Dimensional Crossover and Results}
\hspace{0.5cm}
A major problem in percolation theory is the existence and determination of critical exponents. For instance,  quantities such as $\chi_p(d)$ are believed to diverge as $p\uparrow p_c(d)$ in the manner of a power law in $|p-p_c(d)|$, whose exponent is called a critical exponent (see Chapter 9 in \cite{GRI} for details). More precisely, it is believed that there exists a $\gamma=\gamma(d)>0$ such that 
\begin{equation*}\label{crit}
\chi_p(d)\approx |p-p_c(d)|^{-\gamma},
\end{equation*}
when $p\uparrow p_c(d)$. Here the relation $a(p) \approx b(p)$ means $\log$ equivalence, i.e., $\frac{\log a(p)}{\log b(p)}\rightarrow 1$ when ${p\uparrow p_c(d)}$.

In \cite{LS} the authors introduce another critical exponent, the so-called dimensional crossover critical exponent for the Ising Model, which is related to the function in Equation (\ref{qcrit}). The same exponent is introduced in \cite{RS} for bond percolation and it is expected that:

\begin{Conjecture}\label{crit2} There exists a critical exponent $\psi=\psi(d)> 0$, depending only on $d$, such that
$$q_c(p)\approx |p-p_c(d)|^{\psi}.$$ Moreover, if $\gamma(d)$ exists, then $\psi(d)=\gamma(d)$.
\end{Conjecture}

We stress that the crossover critical exponent $\psi$ should depend only on $d$ and not on $s$. 

We highlight a few papers that investigate this matter.  In \cite{GCGR} the authors examine bond percolation on $\Z^3=\Z^2 \times \Z$ . Here $\Z^2$-edges are open with probability $p$ and $\Z$-edges are open with probability $q=Rp$, where $R$ is the anisotropy parameter. By means of a simulation, the authors estimate $\psi(2)$ by $2.3 \pm 0.1$, which is compatible with the critical exponent $\gamma(2)$, which is expected to be $\frac{43}{18}$ (see \cite{S} for example).
In \cite{RS} the authors study a percolation process on $\Z^d=\Z^{d-1}\times\Z$ where $\Z^{d-1}$-edges are open with probability $p$ and $\Z$-edges parallel to $z$ are open with probability $q= Rp$.  Simulated data then indicate that in the limit  $1/R\rightarrow 0$,  the crossover exponent $\psi$ is equal to $1$ for all $d$. In the opposite limit $R\rightarrow 0$, their analysis suggests that $\psi(d-1)\neq\gamma(d-1)$. This result was later contradicted by Redner and Coniglio \cite{RC}, where the authors argue the opposite relation, that is,  $\psi(d-1)=\gamma(d-1)$.

In \cite{SS,cSS} the authors proved that, if $\gamma(d)$ and $\psi(d)$ exist, then  $\psi(d)\leq \gamma(d)$. In this paper we are concerned with the reversed inequality. The following theorem gives an upper bound for the critical curve $q_c(p)$ when $p$ is sufficiently close to $p_c(d)$, providing a partial answer in that direction.

\begin{theo}\label{teo} Consider an anisotropic bond percolation process on $\Z^d\times\Z^s$, $d,s\geq1$, with parameters $(p,q)$ and $p_c(d)-p>0$,  sufficiently small. If the pair $(p,q)$ satisfies $$q\geq 8d^2(p_c(d)-p),$$ then there is a.s. an infinite open cluster in $\Z^{d+s}$.
\end{theo}

 Theorem \ref{teo} gives an upper bound for $q_c(p)$, \textit{i.e.}, $q_c(p)\leq  8d^2(p_c(d)-p).$ This in turn gives 
 $$\frac{\log(q_c(p))}{\log(p_c(d)-p)}\geq \frac{\log(8d^2)}{\log(p_c(d)-p)}+1. $$ Taking  limits when $p\uparrow p_c(\Z^d)$ we obtain
\begin{equation}\label{liminf}\liminf_{p\uparrow p_c(d)}\frac{\log(q_c(p))}{\log|p-p_c(d)|}\geq 1.\end{equation} In \cite{SS,cSS} it is shown that $$\limsup_{p\uparrow p_c(d)}\frac{\log(q_c(p))}{\log|p-p_c(d)|}\leq \gamma(d),$$ whenever $\gamma(d)$ exists.

The results in \cite{HS} imply that $\gamma(d)$ exists and equals 1 for sufficiently high dimensions, (the authors explicitly prove that $d\geq 19$ is enough in \cite{HS2}). In a more recent work (see \cite{FH}), the same result was shown to hold for $d\geq 11$.  As a consequence, it follows that $\psi(d)$ exists and is equal to 
$\gamma(d)$ for all $d\geq 11$. We have just proved the following: 

\begin{theo}\label{teo2} Consider an anisotropic bond percolation process on $\Z^d\times\Z^s$, $d\geq 11, s\geq1$, with parameters $(p,q)$. Then the critical exponent 
$\psi(d)$ exists and is equal to $\gamma(d)=1$. Hence Conjecture  \ref{crit2} is true in case $d\geq 11$.
\end{theo}

\begin{Rem}Using Theorem 1 of \cite{SS} and the first inequality of $(1.5)$ in \cite{AN}, together with  Theorem \ref{teo} above, we obtain that the ratio $q_c(p)/|p-p_c|$ is bounded and positive, which is stronger than log-equivalence.

\end{Rem}

For spread-out models it is known that the mean-field estimates hold for $d> 6$ (see \cite{HS}) and our results extend to this case. We discuss this in more detail in the last section.

\section{Proof of main result }

The idea of the proof is to construct a dynamical coupling of two percolation processes. This type of construction goes back to the 80's and 90's, see for instance \cite{BS}, \cite{CCN} and \cite{GM}. The proof takes into consideration the exploration of the cluster of the origin of $\Z^{d+s}$ and compare it to a supercritical homogeneous  cluster on $\Z^d$. We see $\Z^{d+s}$ as an ordered layered graph and explore it in the $\Z^d$ direction. Each time we find a closed edge we shift layers with the aid of $\Z^s$-edges. By tuning the parameters $p$ and $q$ properly we show that the percolation process in $\Z^{d+s}$ stochastically dominates supercritical percolation on $\Z^d$, which will give the desired result. Let us dive into the details. 

\textit{Proof of Theorem 1.}
It is sufficient to prove the theorem for the case $s=1$. Let $U$ be the set of unit vectors in $\Z^d$, so that $|U|=2d$. We will denote the vertices of $\Z^{d+1}$ by $(\u,t)$ where $\u$ and $t$ are  the $\Z^d$ and $\Z$ component respectively.



Consider now  the multigraph obtained from the vertices of $\Z^{d+1}=\Z^d \times \Z$ where every $\Z$-edge is replaced by $2d$ other edges indexed by $U$. We denote this graph by $\Z^{d+1}_U$ and introduce percolation on it as follows: as before, every $\Z^d$-edge is open with probability $p$ and every $\Z$-edge is open with probability $\bar{q}$ independently of all others, with $\bar{q}$ satisfying
\begin{equation}\label{alfa}
(1-q)=(1-\bar{q})^{2d}.
\end{equation} 

It is clear that with these parameters, the distribution of the cluster of the origin in $\Z^{d+1}$ with law $\PP_{p,q}$ is the same as that in $\Z^{d+1}_U$ with law $\PP_{p,\bar{q}}$.


As mentioned before, the proof consists  of a dynamical coupling between two percolation processes. That is, for every configuration $\omega\in \{0,1\}^{E(\Z^{d+1}_U)}$ with law $\PP_{p,\bar{q}}$, we shall obtain a configuration $\omega'\in \{0,1\}^{E(\Z^{d})}$ with law $\PP_{r}$ on $\Z^d$ where $r=p+\bar{q} p(1-p)$. To accomplish this we construct  a sequence $\{\eta(e)\}_{ {e}\in E(\Z^d)}$ of independent $0$-$1$ valued random variables with parameter $r$, a sequence  $E_i=(A_i,B_i)$ of ordered pairs of subsets of $E(\Z^d)$, a sequence $S_i$ of subsets of $\Z^{d+1}$ and a sequence $S_i^\pi$ of subsets of $\Z^d$, with $i\in \N$.  We proceed as follows.

Given $(\u,t)\in\Z^{d+1}_U$, we use the notation $\langle(\u,t),(\u,t+1)\rangle_v$ for the $\Z$-edge between $(\u,t)$ and $(\u,t+1)$ indexed by $v\in U$. Recall that $U$ is the set of unit vectors in $\Z^d$.  We say  there is a \textit{$v$-hook} at vertex $(\u,t)\in\Z^{d+1}$ if the edges $\langle(\u,t),(\u,t+1)\rangle_v$ and $\langle(\u,t+1),(\u+v,t+1)\rangle$ are open.
We consider an arbitrary, but fixed, ordering of  $E(\mathbb{Z}^{d})$ and let $(\mathbf{0},0)$ be the origin of $\mathbb{Z}^{d+1}$. We set  $E_0=(\emptyset,\emptyset)$,  $S_0=\{(\mathbf{0},0)\}$ and $S_0^{\pi}=\{\mathbf{0}\}$. The sets $S_k$ and $S^{\pi}_k$ to be defined will be the explored vertices in the open cluster of the origin of $\Z^{d+s}$ and its projection on $\Z^d$, respectively. Let $f_1=\langle\mathbf{0},\mathbf{0}+v\rangle$, $v\in U$, be the first $\Z^d$-edge, in the fixed ordering, incident to $S_0^{\pi}$. We write $[{\bf e},t]\in E(\Z^{d+1})$ for the edge $\langle(\u_1,t),(\u_2,t)\rangle$ whenever ${\bf e}=\langle\u_1,\u_2\rangle\in E(\Z^d)$. Set $\eta(f_1)=1$ if exactly one of the following two conditions hold:

\begin{enumerate}[(a)]
\item $[f_1,0]$ is open;
\item $[f_1,0]$ is closed and there is a \textit{$v$-hook} at vertex $(\mathbf{0},0)$. 
\end{enumerate}
We set
$$E_1=\left\{\begin{array}{ll}
(f_1,\emptyset)&\mbox{if}\quad \eta(f_1)=1,\\
(\emptyset,f_1) &\mbox{if}\quad \eta(f_1)=0,
\end{array}\right.
$$
 If $\eta(f_1)=1$, we set
$$S_1=\left\{\begin{array}{ll}
S_0\cup\{(v,0)\}&\mbox{if  condition (a) holds},\\
S_0\cup\{(v,1)\}&\mbox{if condition (b) holds},
\end{array}\right.
$$
and 
$$S_1^{\pi}=S_0^{\pi}\cup\{v\}.$$
Whenever $\eta(f_1)=0$, we set  $S_1=S_0$ and $S_1^{\pi}=S_0^{\pi}$.

Suppose the sequences $\{E_i\},\{S_i\}$ and $\{S^{\pi}_i\}$ are defined up to the index $i=n-1$. We then define $E_{n}$, $S_n$ and $S^\pi_{n}$ as follows. At first, let $\pi:\Z^{d}\times\Z\rightarrow \Z$ be the projection of $\Z^{d}\times\Z$ onto $\Z$, that is,  $\pi((\u,t))=t$, and consider the function between $S_{n-1}$ and $S^{\pi}_{n-1}$ given by
\begin{align*}
      \tau \colon &S_{n-1} \longrightarrow S_{n-1}^{\pi},\\
       &(\mathbf{u},t) \xmapsto{\phantom{L^\infty(T)}} \mathbf{u}.
\end{align*} 

Note that, by construction, $\tau$ is a bijection. Let $f_n$ be the earliest $\Z^d$-edge (in the fixed ordering) at the outer boundary of $S_{n-1}^{\pi}$, \textit{i.e.}, $f_n\cap S_{n-1}^{\pi}\neq \emptyset$, $f_n\cap (S_{n-1}^{\pi})^c\neq \emptyset$ and $f_n\notin A_{n-1}\cup B_{n-1}$.  Assume, with no loss of generality, that $f_n=\langle\u_{n-1},\u_n\rangle$, where $\u_n=\u_{n-1}+v$ for some $v\in U$, with $\u_{n-1}\in S_{n-1}^{\pi}$ and $\u_{n-1}+v\in (S_{n-1}^{\pi})^c$.

We set $\eta(f_n)=1$ if exactly one of the following two conditions hold:

\begin{enumerate}[(a)]
\item $[f_n,\pi(\tau^{-1}(\u_{n-1}))]$ is open,
\item $[f_n,\pi(\tau^{-1}(\u_{n-1}))]$ is closed and there is a \textit{v-hook} at vertex $(\u_{n-1},\pi(\tau^{-1}(\u_{n-1})))$. 
\end{enumerate}

We then set
$$E_n=\left\{\begin{array}{ll}
(A_{n-1}\cup f_n,B_{n-1})&\mbox{if}\quad \eta(f_n)=1,\\
(A_{n-1},B_{n-1}\cup f_n) &\mbox{if}\quad \eta(f_n)=0,
\end{array}\right.
$$

If $\eta(f_n)=1$, we set
\begin{equation}\label{Sn}S_n=\left\{\begin{array}{ll}
S_{n-1}\cup\{(\u_n,\pi(\tau^{-1}(\u_{n-1})))\}& \mbox{ if condition (a) holds},\\
S_{n-1}\cup\{(\u_n,\pi(\tau^{-1}(\u_{n-1}))+1)\}&\mbox{ if condition (b) holds},
\end{array}\right.
\end{equation}
and 
$$S_n^{\pi}=S_{n-1}^{\pi}\cup\{\u_n\}.$$

In case $\eta(f_n)=0$, we set $S_n=S_{n-1}$ and $S_n^{\pi}=S_{n-1}^{\pi}$.

Now that the dynamical coupling is well defined, we make some observations.  
By construction, there is a bijection between the sets $A_n$ and $S_n\backslash\{({\bf 0},0)\}$ and the sets $A_n$, $B_n$ and $S_n$ are non-decreasing. So we can define $A_{\infty}:=\lim A_n$, $B_{\infty}:=\lim B_n$ and $S_{\infty}:=\lim S_n.$ We also observe that all edges in $A_\infty$ form a connected set containing the origin of $\Z^d$ and also that $S_\infty$ is a subset of the cluster of the origin of $\Z^{d+1}$ in the process with law $\PP_{p,q}$.

To complete the description of the law on $\{0,1\}^{E(\Z^d)}$, we dispose of a collection of \textit{iid} Bernoulli random variables $\{\eta(e)\}$, for all $e\in E(\Z^d)\backslash (A_\infty\cup B_\infty)$, with parameter $r=p+\bar{q}p(1-p)$, independent from all other random variables used previously. 

Now, since the random variables $\{\eta(e)\}_{e\in E(\Z^d)}$ are independent, the probability measure generated by them is exactly the same as that of an independent bond percolation process on $\Z^d$ with parameter $r=p+\bar{q}p(1-p)$. This means that 
$\PP_{p,q}(|A_{\infty}|=\infty)=\theta(r)$.

Taking $q\geq 8d^2(p_c(d)-p)$, observing that $p<p_c(d)\leq 1/2$ and taking $(p_c(d)-p)$ sufficiently small, say $p\in (\frac{1}{2d},p_c(d))$, we estimate

$$r=p+\bar{q}p(1-p)=p+\left[1-(1-q)^{1/2d}\right]p(1-p)>p+\frac{q}{2d}\frac{1}{4d}\geq p_ c(d).$$

To conclude, we observe that, for these values of $q$, we have $$\theta(p,q)\geq \PP_{p,q}(|S_\infty|=\infty)=\PP_{p,q}(|A_{\infty}|=\infty)=\theta(r)>0.$$
\qed

An interesting feature of this dynamical coupling is that, with a minor modification, the same result holds for the bilayered graph $\Z^d\times \{0,1\}$ rather than the full graph $\Z^{d+1}$. Indeed, in the construction above, we always move one layer up each time we use a {\it v-hook}, but we could simply alternate within two layers. Formally,  the exact same proof holds by writing every sum in the $\mathbb{Z}$ coordinate as a 
''mod 2`` sum.

Our bound $\psi\geq 1$ should be sharp for $d>6$ for the full graph $\Z^d\times\Z$ as well as for layered graphs $\Z^d\times\{0,1,\dots,l\}$.  In low dimensions, however, it is not clear that the crossover exponent for layered graphs is the same as the one for the full graph. It would be interesting to investigate how the isoperimetric profile of the sub-critical clusters affects the exponent for each thickness of a layered graph. In particular, one wonders if the actual critical exponent $\psi$ should be strictly smaller than $\gamma$ for layered graphs for $2\leq d\leq 5$.

\section{Discussion}

Under the hypothesis that $\psi(d)$ exists, the expression in (\ref{liminf}) already shows that $\psi(d)\geq 1$. This bound should saturate when the dimension is above the so-called upper critical dimension for percolation ($d_c=6$), but is not expected to be sharp for $2\leq d\leq 5$.

In \cite{HS}, the authors prove that the mean-field behavior occurs for any dimension $d>6$ as long as we consider some spread-out models that should be in the same universality class as the nearest-neighbor model. In particular, their proof works for finite-range models with sufficiently large range. In this case, our proof of Theorem \ref{teo} works as well, with a minor adaptation. Instead of considering the set of unit vectors, we consider $U$ to be the set of edges emanating from the origin. Since $U$ is finite, we simply replace \eqref{alfa} with $(1-q)=(1-\bar{q})^{|U|}$ and obtain the same result as in Theorem \ref{teo}, albeit with a different constant multiplicative factor. Together with the results in \cite{SS,cSS}, we obtain that the critical exponent $\psi$ exists and equals 1 in this case. The same argument works for finite-range spread-out models on the bilayered graph taking again "mod 2" sums in the $\Z$ coordinate.

It should be mentioned that the mean-field
results in \cite{HS} for spread-out models also apply to some infinite-range models. Our proof does not extend to this case.

We also observe that the bound obtained in Theorem 1 is not directly related to $\chi_p(d)$ and we expect that the following result should be true for all $d$.

\begin{Conjecture}\label{teo1} Consider an anisotropic bond percolation process on $\Z^d\times\Z^s$, $d,s\geq1$, with parameters $(p,q)$ and $p_c(d)-p>0$, sufficiently small. There exists a constant $\beta$ such that if the pair $(p,q)$ satisfies $$q>\frac{\beta}{\chi_p(d)},$$ then there is a.s. an infinite open cluster in $\Z^{d+s}$.
\end{Conjecture}

In mean-field, that is, when we consider regular trees instead of $\Z^d$, it is straightforward to prove that the above conjecture holds. We decided not to write down the proof in this note since it only works for trees and does not shed any light in the case of $\Z^d$.

\section{Conflict of Interest}

On behalf of all authors, the corresponding author states that there is no conflict of interest.

\section*{Acknowledgements} The authors would like to thank two anonymous referees for their valuable comments. Remy Sanchis was partially supported by Conselho Nacional de Desenvolvimento Cient\'ifico e Tecnol\'ogico (CNPq) and by  Funda\c c\~ao de Amparo \`a Pesquisa do Estado de Minas
Gerais (FAPEMIG), grant PPM 00600/16. This study was financed in part by the Coordenação de Aperfeiçoamento de Pessoal de Nível Superior – Brasil (CAPES) – Finance Code 001.

\end{document}